\documentclass[9pt]{article}
\usepackage{amssymb,enumerate}
\usepackage{amsmath}
\usepackage{graphics}
\usepackage{cite}
\newtheorem{definition}{\textbf{Definition}}[section]
\newtheorem{theorem}{\textbf{Theorem}}[section]
\newtheorem{lemma}{\textbf{Lemma}}[section]

\begin{document}

\title{Symmetry Analysis of Initial and Boundary Value Problems for Fractional Differential Equations in Caputo sense}
\author{Dogan Kaya$^1$, Gulistan Iskenderoglu$^2$}

\maketitle
\pagestyle{myheadings}
\markboth{Dogan Kaya, Gulistan Iskenderoglu}{Symmetry Analysis of IBVP for Fractional Differential Equations in Caputo sense }


\begin{abstract}
In this work we study Lie symmetry analysis of initial and boundary value problems for partial differential equations (PDE) with Caputo fractional derivative. We give generalized definition and theorem for the symmetry method for PDE with Caputo fractional derivative, according to Bluman's definition and theorem for the symmetry analysis of PDE system. We investigate the symmetry analysis of initial and boundary value problem for fractional diffusion and the third order fractional partial differential equation (FPDE). Also we give some solutions.
\end{abstract}

\section{Introduction}
Theory of differential equations is one of the most important discipline in mathematics that has extensive application in physics.
The descriptions of fundamental physical, chemical, and biological systems phenomena and processes in nature are made by using nonlinear equations. Following to Nobel Laureate Werner Heisenburg's words "... the progress of physics will to a large extent depend on the progress of nonlinear mathematics, of methods to solve nonlinear equations ... and therefore we can learn by comparing different nonlinear problems", - we  will concentrate our attention on nonlinear PDE with fractional derivative \cite{wilhems}.
Numerous complex nonlinear phenomena and dynamic processes in nature can be described or modeled by dint of a FPDE or a system of FPDE which have been used successfully in physics, electromagnetics \cite{zel_mil}, acoustics, astrophysics \cite{taras}, viscoelasticity, chemistry, electrochemistry, etc. \cite{He,Diethelm,MillerRoss,wulee,BarlTorv}. And there are numerous methods to find the solutions for the differential equations like: Adomian decomposition method \cite{Adom,Kaya}, differential transform method \cite{HereBanCha}, modified simple equation method \cite{TagMirRahAkb}, Lie symmetry analysis \cite{HuShZh,SahBak} and so on.

In this work we focus on initial and boundary value problems for FPDE with fractional derivative in the sense of Caputo \cite{Podl}, that has a form
\begin{gather*}
^CD^\alpha_t f(t)=\left\{
\begin{array}{l}
\frac{1}{\Gamma(n-\alpha)}\int_a^t (t-\tau)^{n-\alpha-1} f^{(n)}(\tau)d\tau, \quad\quad if\quad n-1<\alpha<n, n\in\mathbb{N},\\
\mbox{}\\
\frac{d^n}{dt^n}f(t),\quad\quad\quad\quad\quad\quad\quad\quad\quad\quad\quad\quad\quad if\quad\alpha=n.\\
\end{array}
\right.
\end{gather*}
The FPDE with Caputo derivative is more useful in searching the solution of boundary value problems, than Riemann-Liouville derivative, which has a form
\begin{gather*}
D^\alpha_t f(t)=\left\{
\begin{array}{l}
\frac{1}{\Gamma(n-\alpha)}\frac{d^n}{dt^n}\int_a^t (t-\tau)^{n-\alpha-1} f(\tau)d\tau, \quad\quad if\quad n-1<\alpha<n, n\in\mathbb{N},\\
\mbox{}\\
\frac{d^n}{dt^n}f(t),\quad\quad\quad\quad\quad\quad\quad\quad\quad\quad\quad\quad\quad if\quad\alpha=n.\\
\end{array}
\right.
\end{gather*}
 We can assume that Caputo derivative is more handy since the initial value for fractional differential equation with Caputo derivative is the same as the initial value for integer PDE \cite{li_qi_che}. For instance, the initial value condition for fractional differential equation $^CD_{t}^{\alpha}f(t)=u(t,f)$ with $\alpha \in (0,1)$, $t < 0$ is given as $f(0)=f_0$. Along with it for the fractional differential equation $D_{t}^{\alpha}f(t)=u(t,f)$ with $\alpha \in (0,1)$, $t < 0$, the initial value condition can contain fractional integral or derivative and it is given as $[D_{t}^{\alpha-1}f(t)]_{t=0}=f_0$ if $\alpha \in (1,2)$, then the initial value conditions are given as  $[D_{t}^{\alpha-2}f(t)]_{t=0}=f_0$, $[D_{t}^{\alpha-1}f(t)]_{t=0}=f_0^{(1)}$. So, we consider the boundary value problems with Caputo fractional derivative.

It is well-known that Lie group method is a powerful and direct approach to construct exact solutions of nonlinear PDE, by analysing the symmetries of the nonlinear PDE. In addition, based on  Lie symmetry method, many other types of exact solutions of PDE can be obtained, such as traveling wave solutions, soliton solutions, power series solutions, and so on \cite{GKL,KhaAdem,iskndrkaya}. Lie symmetry method was investigated by Sophus Lie, later Ovsyannikov, Olver and many other big minds studied it by contributing to the theory of symmetries \cite{Ovsyannikov,Ibrag1,olver}. Bluman investigated the symmetry analysis of initial and boundary value problems for PDE, by given the properties of boundary conditions and boundary surfaces under symmetries (see \cite{BlumKumei,BlumAnco}, and references therein). And we in our work by using Lie symmetry method and its analysis search symmetries and obtain solutions for initial and boundary value problems for PDE, containing Caputo fractional derivative.

\section{Symmetry analysis for fractional partial differential equations}

Consider a time-fractional PDE with two independent variables $x \in \mathbb{R}$ and $t>0$ is given as following:
\begin{equation}\label{fpde}
F(x, t, u, {^CD}^\alpha_t u, \partial_x u, \partial_x^2 u,...,\partial_x^s u)=0, \quad  0<\alpha\leq1
\end{equation}
where ${^CD}^\alpha_t u$ is Caputo fractional derivative of $u$ and $\partial_x^i u=\frac{\partial^i u}{\partial x^i}$, $i=1,...,s$.

One point Lie symmetry transformation acting on a space of two independent variables $(t,x)$ and depended variable $u$ is determined as
\begin{equation}\label{Lietrans}
\begin{array}{l}
\bar{t}=t+\varepsilon \tau(t,x,u)+O(\varepsilon^2),\\
\mbox{}\\
\bar{x}=x+\varepsilon \xi(t,x,u)+O(\varepsilon^2),\\
\mbox{}\\
\bar{u}=u+\varepsilon \eta(t,x,u)+O(\varepsilon^2),\\
\end{array}
\end{equation}
where $\varepsilon > 0$ is an infinitesimal group parameter and $\tau, \xi, \eta$ are the infinitesimals of the transformation. An infinitesimal generator associated with the above transformation has a following form
\begin{equation}\label{GENER}
X=\xi(t,x,u)\frac{\partial}{\partial x}+\tau(t,x,u)\frac{\partial}{\partial t}+\eta(t,x,u)\frac{\partial}{\partial u}
\end{equation}
where $\xi=\frac{d\bar{x}}{{d\varepsilon}}\mid_{\varepsilon=0}$, $\tau=\frac{d\bar{t}}{{d\varepsilon}}\mid_{\varepsilon=0}$ and $\eta=\frac{d\bar{u}}{{d\varepsilon}}\mid_{\varepsilon=0}$.
\begin{definition} We say $u=\theta(t,x)$ is an invariant solution of the equation (\ref{fpde}) resulting from the invariance of the equation (\ref{fpde}) under the symmetry (\ref{Lietrans}) with infinitesimal generator (\ref{GENER}) if only if (see \cite{olver})
\begin{itemize}
\item $u=\theta(t,x)$ is an invariance surface of $X$,
\item $u=\theta(t,x)$ solves the equation (\ref{fpde}).
\end{itemize}
\end{definition}
In other words, we can say that $u=\theta(t,x)$ is the invariant solution of the equation (\ref{fpde}) under the infinitesimal generator $X$ with Lie transformation (\ref{Lietrans}) if only if $u=\theta(t,x)$ satisfies:
\begin{itemize}
\item $X(u-\theta(t,x))=0$ for $u=\theta(t,x)$, which gives us
$$
\xi(t,x,\theta(t,x))\frac{\partial \theta(t,x)}{\partial x}+\tau(t,x,\theta(t,x))\frac{\partial \theta(t,x)}{\partial t}=\eta(t,x,\theta(t,x)),
$$
\item $F(x, t, \theta(t,x), {^CD}^\alpha_t \theta(t,x), \partial_x \theta(t,x), \partial_x^2 \theta(t,x),...,\partial_x^s \theta(t,x))=0$, for $u=\theta(t,x)$.
\end{itemize}
According to the infinitesimal transformation (\ref{Lietrans}) a fractional prolongation $pr^{(\alpha, n)}X$ of the equation (\ref{fpde})
$$
pr^{(\alpha, s)}X(E)\mid_{E=0}=0, \quad \quad where\quad  E=F(x, t, u, \partial^\alpha_t u, \partial_x u, \partial_x^2 u,...,\partial_x^s u)=0,
$$
is an invariant criteria for equation (\ref{fpde}) and has the form
$$
pr^{(\alpha, n)}X=X+\eta^t_\alpha\partial_{{^CD}^\alpha_t u}+\eta^x_1\partial_{u_x}+\eta^x_2\partial_{\partial_x^2 u}+...+\eta^x_s\partial_{\partial_x^s u},
$$
here
\begin{equation}\label{eta}
\begin{array}{c}
\eta^t_\alpha={^CD}^{\alpha}_t(\eta)+\xi {^CD}^{\alpha}_t(u_x)-{^CD}^{\alpha}_t(\xi u_x)+\tau ^CD^{\alpha}_t(u_t)-{^CD}^{\alpha}_t(\tau u_t),\\
\mbox{}\\
\eta^x_1=D_x\eta-u_xD_x\xi_1-u_tD_x\tau,\\
\mbox{}\\
\eta^x_2=D_x\eta^x_1-u_{xx}D_x\xi-u_{xt}D_x\tau,\\
\mbox{}
\vdots
\mbox{}\\
\eta^x_s=D_x\eta^x_{s-1}-\partial_x^s u D_x\xi-\partial_x^{s-1} u_t D_x\tau,\\
\end{array}
\end{equation}
with total derivative $D_i$ in a form
$$
D_i=\partial_i+u_i\partial_u+u_{it}\partial_{u_t}+u_{jt}\partial_{u_t}+u_{ii}\partial_{u_i}+u_{jj}\partial_{u_j}+... .
$$

The expressions for $\eta^x_i$, $i=1,...,s$ in (\ref{eta}) can be easily obtained (see \cite{olver}), here we focus our attention on $\eta^t_\alpha$ \cite{olver,Ibrag} .

We can define the relation between Caputo fractional derivative and Riemann--Liouville fractional derivative in a form
$$
D^\alpha_t f(t) =^CD^\alpha_t f(t)+\sum\limits_{i=0}^{n-1}\frac{f^{(i)}(0)}{\Gamma(i-\alpha+1)}t^{i-\alpha}.
$$
And the fractional integral can be defined as
$$
I^\alpha_t f(t)=\frac{1}{\Gamma(\alpha)}\int_0^t (t-\tau)^{\alpha-1}f(\tau)d\tau.
$$
So, we can see the below representations of Caputo and Riemann--Liouville fractional derivatives (see \cite{MillerRoss,Podl})
\begin{equation}\label{CI}
{^CD}^\alpha_tf(t)=I^{n-\alpha}D^nf(t)
\end{equation}
\begin{equation}\label{RLI}
{D}^\alpha_tf(t)=D^nI^{n-\alpha}f(t)
\end{equation}

According to the generalized Leibnitz rule in \cite{SamKilMar}
$$
D^\alpha_t(f(t)g(t))=\sum\limits^{\infty}_{n=0}
\left(\begin{matrix}
\alpha\\
n
\end{matrix}
\right)
I^{n-\alpha}_tf(t)D^n_tg(t), \quad \quad \left(\begin{matrix}
\alpha\\
n
\end{matrix}
\right) =\frac{(-1)^{n-1}\alpha\Gamma(n-\alpha)}{\Gamma(1-\alpha)\Gamma(n+1)},
$$
we have
$$
\xi {^CD}^{\alpha}_t(u_x)-{^CD}^{\alpha}_t(\xi u_x)=-\sum\limits^{\infty}_{n=1}
\left(\begin{matrix}
\alpha\\
n
\end{matrix}
\right)
I^{n-\alpha}_t(u_x)D^n_t(\xi)+\sum\limits_{i=0}^{n-1}\frac{(\xi u_x)^{(i)}(0)}{\Gamma(i-\alpha+1)}t^{i-\alpha},
$$
and
$$
\tau{^CD}^{\alpha}_t(u_t)-{^CD}^{\alpha}_t(\tau u_t)=-\sum\limits^{\infty}_{n=1}
\left(\begin{matrix}
\alpha\\
n
\end{matrix}
\right)
I^{n-\alpha}_t(u_t)D^{n}_t(\tau)+\sum\limits_{i=0}^{n-1}\frac{(\tau u_t)^{(i)}(0)}{\Gamma(i-\alpha+1)}t^{i-\alpha}.
$$
Thus, we get the expression
$$
\begin{array}{l}\eta^t_\alpha=D^{\alpha}_t(\eta)+\eta_0-\sum\limits^{\infty}_{n=1}
\left(\begin{matrix}
\alpha\\
n
\end{matrix}
\right)
I^{n-\alpha}_t(u_x)D^n_t(\xi)-\sum\limits_{i=0}^{n-1}\frac{(\xi u_x+\tau u_t)^{(i)}(0)}{\Gamma(i-\alpha+1)}t^{i-\alpha}-\\
\sum\limits^{\infty}_{n=1}
\left(\begin{matrix}
\alpha\\
n+1
\end{matrix}
\right)
I^{n-\alpha}_t(u)D^{n+1}_t(\tau)-\alpha D_t(\tau)D_t^\alpha(u).
\end{array}
$$
Here by using a generalized form of the chain rule (see \cite{Oldham})
$$
\frac{d^mf(g(t))}{dt^m}=\sum\limits^{m}_{k=0}\sum\limits^{k}_{r=0}
\left(\begin{matrix}
k\\
r
\end{matrix}
\right)
\frac{1}{k!}\left(-g(t)\right)^r\frac{d^m}{dt^m}\left(g(t)^{k-r}\right)\frac{d^kf(g)}{dg^k},
$$
which is
{\footnotesize$$
\frac{d^\alpha f(t, g(t))}{dt^\alpha}=\sum\limits^{\infty}_{l=0}\sum\limits^{l}_{m=0}\sum\limits^{m}_{k=0}\sum\limits^{k-1}_{r=0}
\left(\begin{matrix}\alpha\\l\end{matrix}\right)\left(\begin{matrix}l\\m\end{matrix}\right)\left(\begin{matrix}k\\r\end{matrix}\right)
\frac{1}{k!}\frac{t^{l-\alpha}}{\Gamma(l+1-\alpha)}\left(-g(t)\right)^r\frac{d^m}{dt^m}\left(g(t)^{k-r}\right)\frac{\partial^{l-m+k}f(g)}{\partial t^{l-m}\partial g^k},
$$}
we have
$$
\eta_0=\sum\limits_{j=0}^{n-1}\frac{t^{n-\alpha}(\eta)^{(j)}(0)}{\Gamma(j+1-\alpha)}= \sum\limits_{j=0}^{n-1}\frac{t^{(n-\alpha)}}{\Gamma(j+1-\alpha)}\sum\limits^{\infty}_{l=0}\sum\limits^{l}_{m=0}\sum\limits^{m}_{k=0} \sum\limits^{k-1}_{r=0}
\left(\begin{matrix}\alpha\\l\end{matrix}\right)\left(\begin{matrix}l\\m\end{matrix}\right)\left(\begin{matrix}k\\r\end{matrix}\right)
\frac{1}{k!}\times
$$
$$
\frac{t^{l-\alpha}}{\Gamma(l+1-\alpha)}\left(\left(-u\right)^r\frac{d^j}{dt^j}\left(u^{k-r}\right)\frac{\partial^{l-j+k}\eta}{\partial t^{l-j}\partial u^k}\right)(0).
$$
And so, the infinitesimal $\eta^t_\alpha$ takes a form
\begin{equation}\label{etat1}
\begin{array}{l}
\eta^t_\alpha=\frac{^C\partial^\alpha \eta}{\partial t^\alpha} + (\eta_u-\alpha(\tau_t+u_t\tau_u))\frac{^C\partial^\alpha u}{\partial t^\alpha}-u\frac{^C\partial^\alpha \eta_u}{\partial t^\alpha} + \mu-\sum\limits^{\infty}_{n=1}
\left(\begin{matrix}
\alpha\\
n
\end{matrix}
\right)
D^n_t(\xi)I^{n-\alpha}_t(u_x)+\\\sum\limits^{\infty}_{n=1}\left[
\left(\begin{matrix}
\alpha\\
n
\end{matrix}
\right)
\frac{\partial^n \eta_u}{\partial t^n}-
\left(\begin{matrix}
\alpha\\
n+1
\end{matrix}
\right)
D^{n+1}_t\tau\right]I^{n-\alpha}_tu+\sum\limits_{i=0}^{n-1}\frac{(\xi u_x)^{(i)}(0)}{\Gamma(i-\alpha+1)}t^{i-\alpha},
\\
\end{array}
\end{equation}
where
\begin{equation}\label{mu}
\mu=\sum\limits^{\infty}_{n=2}\sum\limits^{n}_{m=2}\sum\limits^{m}_{k=2}\sum\limits^{k-1}_{r=0}\left(\begin{matrix}
\alpha\\
n
\end{matrix}
\right)\left(\begin{matrix}
n\\
m
\end{matrix}
\right)\left(\begin{matrix}
k\\
r
\end{matrix}
\right)\frac{1}{k!}\frac{t^{n-\alpha}(-u)^r}{\Gamma(n+1-\alpha)}\frac{\partial^m}{\partial t^m}\left(u^{k-r}\right)\frac{\partial^{n-m+k}\eta}{\partial t^{n-m}\partial u^k}.
\end{equation}

Here we can proof the next lemma.
\begin{lemma}
If (\ref{fpde}) is invariant under infinitesimal transformation (\ref{Lietrans}) with infinitesimal generator (\ref{GENER}) and the equation (\ref{fpde}) has no the second and higher order derivative of $u$ with respect to $t$, then $\eta=A(t,x)u+B(t,x)$ with $A(t,x)$ and $B(t,x)$ arbitrary functions.
\end{lemma}
\textbf{Proof.} By expanding the ${^CD}^\alpha_t\eta$ we get the expression (\ref{mu}). As the equation has not any variations of second and high order derivative $u$ with respect to $t$, then for $k=2$ we have
$$
\mu=\frac{1}{2!}\frac{t^{2-\alpha}(-u)}{\Gamma(3-\alpha)}u_{tt}\eta_{uu},
$$
and $\eta_{uu}=0$. Here by integration we get $\eta=A(t,x)u+B(t,x)$.
\\Q.E.D.

\section{Symmetry analysis for time-fractional initial and boundary value problems}

In this section we study the one point Lie symmetry analysis for initial and boundary value problems. We know that the PDE can describe real processes according nature and society if there are given initial and boundary conditions for the PDE.
Although Lie symmetry analysis is one of the most widely-applicable methods of finding exact solutions of differential equations, but it was not widely used for solving boundary value problems. The reason is the initial and boundary conditions usually are not invariant under any obtained Lie symmetry method transformations \cite{BlumKumei, BlumAnco}.
Thereby, an invariant solution for PDE resulting by applying symmetry transformation solves a given boundary value problem, when the symmetry transformation leaves invariant all boundary conditions and the domain of the boundary value problem \cite{BlumKumei}.

In \cite{BlumKumei} Bluman gives a definition of Lie symmetry invariance for initial and boundary by means of that there are some classes initial and boundary problems which can be solved. Bluman studies the PDE with finite values of $x$ and $t$ boundary and initial conditions. Later Cherniha extends the definitions of Bluman for the PDE with infinite values of $x$ and $t$ boundary and initial conditions \cite{CherKov, CherKov2}. In this work we investigate the symmetry analysis of initial and boundary value problems for fractional nonlinear PDE with Caputo time-fractional derivative.

Let Lie symmetry infinitesimal generator $X$
\begin{equation}\label{X}
X=\xi(x,t,u) \frac{\partial}{\partial x}+\tau(x,t,u)\frac{\partial}{\partial t}+\eta(x,t,u)\frac{\partial}{\partial u}
\end{equation}
is admitted by the boundary value problem defined on a domain $\Omega$:

\begin{equation}\label{BVPoblem}
u_t=f\left(x,u,\frac{\partial u}{\partial x},...,\frac{\partial^{k} u}{\partial x^{k}}\right), \quad \quad (x,t) \in\Omega \subset \mathbb{R}^2,
\end{equation}
\begin{equation}\label{BV}
d_a(x,t)=0 : B_a\left(t,x,u,\frac{\partial u}{\partial x},...,\frac{\partial^{k-1} u}{\partial x^{k-1}}\right)=0, \quad \quad a=1,...,p.
\end{equation}

Here $B_a(t,x,u,\frac{\partial u}{\partial x},...,\frac{\partial^{k-1} u}{\partial x^{k-1}})$ boundary condition on $d_a(x,t)$. Suppose that the above boundary value problem has a unique solution.

\begin{definition} \cite{BlumKumei} The symmetry $X$ which has the form (\ref{X}) is allowed by the boundary value problem (\ref{BVPoblem}, \ref{BV}) if:
\begin{itemize}
  \item $X^{(k)}(u_t-f(x,u,\frac{\partial u}{\partial x},...,\frac{\partial^{k} u}{\partial x^{k}}))=0$ for $u_t=f(x,u,\frac{\partial u}{\partial x},...,\frac{\partial^{k} u}{\partial x^{k}})$;
  \item $X d_a(x,t)=0$ for $d_a(x,t)=0$, $a=1,...,p$;
  \item $X^{(k-1)}B_a(t,x,u,\frac{\partial u}{\partial x},...,\frac{\partial^{k-1} u}{\partial x^{k-1}})=0$ for $B_a(t,x,u,\frac{\partial u}{\partial x},...,\frac{\partial^{k-1} u}{\partial x^{k-1}})=0$ on $d_a(x,t)=0$, $a=1,..,p$.
\end{itemize}
\end{definition}

Further we extend the Blusman's definition for FPDE with Caputo derivative. Let us consider the boundary value problem for FPDE  defined on a domain $\Omega$

\begin{equation}\label{FBVPoblem}
^CD^\alpha_tu=g\left(x,u,\frac{\partial u}{\partial x},...,\frac{\partial^{k} u}{\partial x^{k}}\right), \quad \quad (x,t) \in\Omega \subset \mathbb{R}^2,
\end{equation}
\begin{equation}\label{FBV}
{^Cd}_a(x,t)=0 : {^CB}_a\left(t,x,u,\frac{\partial u}{\partial x},...,\frac{\partial^{k-1} u}{\partial x^{k-1}}\right)=0, \quad \quad a=1,...,p.
\end{equation}

Here ${^CB}_a(t,x,u,\frac{\partial u}{\partial x},...,\frac{\partial^{k-1} u}{\partial x^{k-1}})$ is a boundary condition on ${^Cd}_a(x,t)$. Suppose the boundary value problem (\ref{FBVPoblem}-\ref{FBV}) has a unique solution. Then we can give next definition.

\begin{definition} The symmetry $X$ in the form (\ref{X}) is allowed by the boundary value problem (\ref{FBVPoblem}, \ref{FBV}) if:
\begin{itemize}
  \item $X^{(k)}({^CD}^\alpha_t u-g(x,u,\frac{\partial u}{\partial x},...,\frac{\partial^{k} u}{\partial x^{k}}))=0$ for ${^CD}^\alpha_t u=g(x,u,\frac{\partial u}{\partial x},...,\frac{\partial^{k} u}{\partial x^{k}})$;
  \item $X {^Cd}_a(x,t)=0$ for ${^Cd}_a(x,t)=0$, $a=1,...,p$;
  \item $X^{(k-1)}{^CB}_a(t,x,u,\frac{\partial u}{\partial x},...,\frac{\partial^{k-1} u}{\partial x^{k-1}})=0$ for ${^CB}_a(t,x,u,\frac{\partial u}{\partial x},...,\frac{\partial^{k-1} u}{\partial x^{k-1}})=0$ on ${^Cd}_a(x,t)=0$, $a=1,..,p$.
\end{itemize}
\end{definition}

\begin{theorem}
The solution $u=\nu(t,x)$ for (\ref{FBVPoblem}) is invariant if only if for infinitesimal generator $X$ the curve (mapping) $\nu(t,x)$ admits
\begin{equation}
\begin{array}{c}
\eta(t,x,\nu(t,x))-\xi(t,x,\nu(t,x))\\
-\tau(t,x,\nu(t,x)){^CD}^{1-\alpha}_tg\left(x,\nu(t,x),\frac{\partial \nu(t,x)}{\partial x},...,\frac{\partial^{k} \nu(t,x)}{\partial x^{k}}\right)=0.
\end{array}
\end{equation}
\end{theorem}
\textbf{Proof:}
As the solution surface  $u=\nu(t,x)$  is invariant for (\ref{FBVPoblem}) if only if  $X(u-\nu(t,x))=0$ which gives
$$0=X(u-\nu(t,x))=\left(\xi(x,t,u) \frac{\partial}{\partial x}+\tau(x,t,u)\frac{\partial}{\partial t}+\eta(x,t,u)\frac{\partial}{\partial u}\right)(u-\nu(t,x))=$$
$$
=\eta(x,t,u)\frac{\partial u}{\partial u}-\xi(x,t,u) \frac{\partial \nu(t,x)}{\partial x}-\tau(x,t,u)\frac{\partial \nu(t,x)}{\partial t}=
$$
$$
=\eta(x,t,\nu(t,x))-\xi(x,t,\nu(t,x)) \frac{\partial \nu(t,x)}{\partial x}-\tau(x,t,\nu(t,x))\frac{\partial \nu(t,x)}{\partial t}.
$$
Further, as $u=\nu(t,x)$, then $u_t=\frac{\partial \nu(t,x)}{\partial t}$. From the property of Caputo derivative (\ref{CI}) $u_t={^CD}_t^{1-\alpha}g(x,u,\frac{\partial u}{\partial x},...,\frac{\partial^{k} u}{\partial x^{k}}))$, i.e
$$
\eta(x,t,\nu(t,x))-\xi(x,t,\nu(t,x)) \frac{\partial \nu(t,x)}{\partial x}-\\
$$
$$
\tau(x,t,\nu(t,x)){^CD}_t^{1-\alpha}g(x,\nu(t,x),\frac{\partial \nu(t,x)}{\partial x},...,\frac{\partial^{k} \nu(t,x)}{\partial x^{k}}))=0.
$$
\\Q.E.D.

\section{Symmetry analysis for some time-fractional initial and boundary value problems}

In this section we will study the symmetries of initial and boundary value problems with fractional derivatives.
Suppose the initial and boundary value problems for fractional diffusion equation
\begin{equation}\label{Diffus}
{^CD}^\alpha_t u=(u^p u_x)_x, \quad\quad\quad p\neq0,
\end{equation}
with initial and  boundary conditions
\begin{equation}\label{IBC}
\left\{\begin{array}{l}
u(t,0)=a(t), \quad\quad for \quad t>0,\\
u(0,x)=b(x), \quad\quad for \quad x>0.
\end{array}\right.
\end{equation}

Gazizov and others in \cite{GKL} found the symmetries for above equation in form

$$
X=(c_1+\frac{\alpha}{2}xc_2+p xc_3)\frac{\partial}{\partial x}+tc_2\frac{\partial}{\partial t}+2uc_3\frac{\partial}{\partial u}.
$$
The invariance of $t=0$ leads to $\tau(0)=0$ and the invariance of $x=0$ leads to $\xi(0,t)=0$ for any $t>0$.
Thus, we get $c_1=0$ for $\xi(0,t)=0$ with $t>0$
$$
tc_2\frac{d}{dt}a(t)=2a(t)c_3,
$$
and for $\tau(0)=0$ we have
$$
\left(\frac{\alpha}{2}xc_2+p xc_3\right)\frac{d}{dx}b(x)=2b(x)c_3,
$$
so, we get
$$
a(t)=t^{\frac{2c_3}{c_2}}k_1, \quad\quad k_1 \quad is \quad a \quad constant,
$$
$$
b(x)=x^{\frac{2}{c_2\alpha/2c_3+p}}k_2, \quad\quad k_2 \quad is\quad a \quad constant.
$$
It means that the initial and boundary value problems can be invariant if the initial and boundary conditions have above forms.

Therefore the equation (\ref{Diffus}) with initial and boundary values (\ref{IBC}) have two symmetry operators
$$
X_1=\frac{\alpha}{2}x\frac{\partial}{\partial x}+t\frac{\partial}{\partial t}, \quad and \quad X_2=p x\frac{\partial}{\partial x}+2u\frac{\partial}{\partial u},
$$
where $X_1$ gives $u=f(x^{\frac{2}{\alpha}}t^{-1})$, and from $X_2$ operator we have $u=x^{\frac{2}{p}}g(t)$. As the problem (\ref{Diffus}) has only one solution, then
$$
f(x^{\frac{2}{\alpha}}t^{-1})=x^{\frac{2}{p}}g(t),
$$
or $f(x^{\frac{2}{\alpha}}t^{-1})=(x^{\frac{2}{\alpha}}g(t)^{\frac{p}{\alpha}})^{\frac{\alpha}{p}}$. Here we see that $f(z)=z^{\frac{\alpha}{p}}$ and $g(t)^{\frac{p}{\alpha}}=t^{-1}$, which gives $g(t)=t^{-\frac{\alpha}{p}}k$. Therefore $u=x^{\frac{2}{p}}t^{-\frac{\alpha}{p}}c$, where $c$ is arbitrary constant, that we will be obtained.

Now let's consider the boundary and initial conditions. Here for $p>0$, $u(t,0)=0$, or in our problem $k_1=0$ and we have blow-up in $t=0$.
Namely for $p>0$, $\lim\limits_{t\rightarrow 0}u(t,x)=\lim\limits_{t\rightarrow 0}x^{\frac{2}{p}}t^{-\frac{\alpha}{p}}c=+\infty$.

And similarly for $p<0$, $u(0,x)=0$, in our problem $k_2=0$ with  blow-up in $x=0$, i.e. $\lim\limits_{x\rightarrow 0}u(t,x)=\lim\limits_{x\rightarrow 0}x^{\frac{2}{p}}t^{-\frac{\alpha}{p}}c=+\infty$.

By checking of the solution for the equation (\ref{Diffus}) we have $k=\sqrt[p]{\frac{2(2-p)\Gamma(1-\frac{\alpha}{p})}{p^2\Gamma(1-\frac{\alpha}{p}-\alpha)}}$, thus
$$
u(t,x)=x^{\frac{2}{p}}t^{-\frac{\alpha}{p}}\sqrt[p]{\frac{2(2-p)\Gamma(1-\frac{\alpha}{p})}{p^2\Gamma(1-\frac{\alpha}{p}-\alpha)}}
$$
is a solution for
$$
\left\{\begin{array}{l}
{^CD}^\alpha_t u=(u^p u_x)_x, \quad\quad\quad p>0,\\
u(t,0)=0, \quad\quad\quad\quad\quad\quad t>0,\\
\lim\limits_{t\rightarrow 0}u(t,x)=+\infty, \quad\quad\quad x>0.
\end{array}\right.
$$
and
$$
\left\{\begin{array}{l}
{^CD}^\alpha_t u=(u^p u_x)_x, \quad\quad\quad p<0,\\
\lim\limits_{x\rightarrow 0}u(t,x)=+\infty, \quad\quad\quad t>0\\
u(0,x)=0, \quad\quad\quad\quad\quad\quad x>0.
\end{array}\right.
$$

Further consider the initial and boundary value problem for below equation

\begin{equation}\label{kdv}
{^CD}^\alpha_t u=(u^q u_{xx})_x, \quad\quad\quad q\neq0,
\end{equation}
with initial and  boundary conditions
\begin{equation}\label{IBC2}
\left\{\begin{array}{l}
u(t,0)=d(t), \quad\quad for \quad t>0,\\
u(0,x)=e(x), \quad\quad for \quad x>0.
\end{array}\right.
\end{equation}

By applying Lie symmetry method for equation (\ref{kdv}) we obtain
\begin{equation}\label{generat}
X=\left(c_1x+c_2\right)\frac{\partial}{\partial x}+\left(c_3t+c_4\right)\frac{\partial}{\partial t}+\left(\frac{3c_1-\alpha c_3}{q}\right)u\frac{\partial}{\partial u},
\end{equation}
Thus the invariance of $t=0$ leads to $c_4=0$ and the invariance of $x=0$ leads to $c_2=0$. And
$$
c_3 t\frac{d}{dt}d(t)=\frac{3c_1-\alpha c_3}{q} d(t),
$$
which gives $d(t)=t^{\frac{3c_1-\alpha c_3}{c_3q}}k_3$, here $k_3$ is an arbitrary constant.

So, for $t=0$ we get
$$
c_1 x\frac{d}{dx}e(x)=\frac{3c_1-\alpha c_3}{q} e(x),
$$
or $e(x)=x^{\frac{3c_1-\alpha c_3}{c_1q}}k_4$, here $k_4$ is an arbitrary constant.

That means the initial and boundary problem for equation (\ref{kdv}) is invariant according to infinitesimal generator (\ref{generat}) if the boundary and initial conditions have above form and there are two infinitesimal operators
$$
Y_1=x\frac{\partial}{\partial x}+\frac{3}{q}u\frac{\partial}{\partial u}, \quad\quad Y_2=t\frac{\partial}{\partial t}-\frac{\alpha}{q}u\frac{\partial}{\partial u}.
$$
The first operator gives $u=t^{-\frac{\alpha}{q}}h(x)$, and the second operator gives $u=x^{\frac{3}{q}}l(t)$. Thus, From the uniqueness of the solution for initial and boundary problem, we have $u=t^{-\frac{\alpha}{q}}x^{\frac{3}{q}}c'$, where $c'$ is an arbitrary constant.

Therefore for $q>0$, $u(t,0)=0$, i.e. $k_3=0$ and we have blow-up in $t=0$ as $\lim\limits_{t\rightarrow 0}u(t,x)=\lim\limits_{t\rightarrow 0}t^{-\frac{\alpha}{q}}x^{\frac{3}{q}}c'=+\infty$.

In a like manner for $q<0$, $u(0,x)=0$, with $k_4=0$ and  blow-up in $x=0$, with $\lim\limits_{x\rightarrow 0}u(t,x)=\lim\limits_{x\rightarrow 0}t^{-\frac{\alpha}{q}}x^{\frac{3}{q}}c'=+\infty$.

By checking the solution for equation (\ref{kdv}) we obtain
$$
 u(t,x)=t^{-\frac{\alpha}{q}}x^{\frac{3}{q}} \sqrt[q]{\frac{q^3\Gamma(1-\frac{\alpha}{q})}{3(3-q)(3+q)\Gamma(1-\frac{\alpha}{q}-\alpha)}},
$$
the solution for initial and boundary value problem
$$
\left\{\begin{array}{l}
{^CD}^\alpha_t u=(u^q u_{xx})_x, \quad\quad\quad q>0,\\
u(t,0)=0, \quad\quad\quad\quad\quad\quad t>0,\\
\lim\limits_{t\rightarrow 0}u(t,x)=+\infty, \quad\quad\quad x>0.
\end{array}\right.
$$
and
$$
\left\{\begin{array}{l}
{^CD}^\alpha_t u=(u^q u_{xx})_x, \quad\quad\quad q<0,\\
\lim\limits_{x\rightarrow 0}u(t,x)=+\infty, \quad\quad\quad t>0,\\
u(0,x)=0, \quad\quad\quad\quad\quad\quad x>0.
\end{array}\right.
$$

\section{Conclusion}
In this work, we present the applications of Lie group analysis to study the initial and boundary value problems for time-fractional nonlinear PDE given by Caputo sense. We give the definition of invariance of the initial and boundary value problems for time-fractional PDE and proved some theorems. By using  Lie symmetry method and its analysis for the initial and boundary value problems with time-fractional PDE we obtained some exact solutions.

\section{Acknowledgment}
The authors would like to thank Istanbul Commerce University, 22-2018/34 and Turkish scholarships for supporting this work.

\end{document}